\documentclass[12pt]{amsart}

\usepackage{amssymb}
\usepackage{xspace}
\usepackage{enumerate}

\newtheorem{thm}{Theorem}

\newtheorem{prop}[thm]{Proposition}

\newtheorem{conj}[thm]{Conjecture}

\theoremstyle{definition}

\theoremstyle{remark}

\newcommand{\Fq} {{\ensuremath{\mathbb{F}_q}}\xspace}

\newcommand{\Fqc} {{\ensuremath{\overline{\,\mathbb{F}\,}_q}}\xspace}
\newcommand{\e} {{\ensuremath{\varepsilon}}\xspace}

\DeclareMathOperator{\diam}{diam}
\DeclareMathOperator{\cl}{cl}

\def\BZ {\mathbb{Z}}

\def\CC {\mathcal{C}}

\begin{document}
\title{Growth in finite simple groups of Lie type} \author{L\'aszl\'o
  Pyber and Endre Szab\'o}
\maketitle

\begin{abstract}
  We prove that if $L$ is a finite simple group of Lie type and $A$ a
  symmetric set of generators of $L$, then $A$ grows i.e $|AAA| > |A|^{1+\e}$
  where $\e$ depends only on the Lie rank of $L$, or $AAA=L$.  This
  implies that for a family of simple groups $L$ of Lie type 
  of bounded rank the
  diameter of any Cayley graph is polylogarithmic in $|L|$.
  \par
  Combining our result on growth with known results of Bourgain,
  Gamburd and Varj\'u it follows that if $\Lambda$ is a Zariski-dense
  subgroup of $SL(d,\BZ)$ generated by a finite symmetric set $S$, then
  for square-free moduli $m$,
  which are relatively prime to some number $m_0$,
  the Cayley graphs $\Gamma(SL(d,\BZ/m\BZ),\pi_m(S))$ form an expander family.
\end{abstract}

\section{Introduction}

The diameter, $\diam(X)$, of an undirected graph $X=(V,E)$ is the
largest distance between two of its vertices.  

Given a subset $A$ of
the vertex set $V$ the expansion of $A$, $c(A)$, is defined to be the
ratio $|\sigma(X)|/|X|$ where $\sigma(X)$ is the set of vertices at
distance $1$ from $A$. A graph is a $C$-expander for some $C>0$ if for
all sets $A$ with $|A|< |V|/2$ we have $c(A)\ge C$. A family of graphs
is an expander family if all of its members are $C$-expanders for some
fixed positive constant $C$.

Let $G$ be a finite group and $S$ a symmetric set of generators of
$G$.  The Cayley graph $\Gamma(G,S)$ is a graph whose vertices are the
elements of $G$ and which has an edge from $x$ to $y$ if and only if
$x=sy$ for some $s \in S$. Then the diameter of $\Gamma$ is the
smallest number $d$ such that $S^d=G$.

The following classical conjecture is due to Babai \cite{BS}

\begin{conj}[Babai]
  For every non-abelian finite simple group $G$ and every symmetric
  generating set $S$ of $G$ we have 
  $\diam\big( \Gamma(G,S)\big) =C\big(\log|G|\big)^c$ 
  where $c$ and $C$ are absolute constants.
\end{conj}

In a spectacular breakthrough Helfgott \cite{He1} proved that the
conjecture holds for the family of groups $G=PSL(2,p)$, $p$ a prime.
In recent major work \cite{He2} he proved the conjecture for the
groups $G=PSL(3,p)$, $p$ a prime.
Dinai \cite{Di} and Varj\'u \cite{Va} have extended
Helfgott's original result to the groups $PSL(2,q)$, $q$ a prime
power.

We prove the following.
\begin{thm}
  Let $L$ be a finite simple group of Lie type of rank $r$. For every
  set $S$ set of generators of $L$ we have 
  $$
  \diam\big(\Gamma(L,S)\big) < C(r)\big(\log|L|\big)^{c(r)}
  $$
  where the constants $C(r)$ and $c(r)$ depend only on $r$.
\end{thm}

This settles Babai's conjecture for any family of simple groups of Lie
type of bounded rank.

\section{Results on growth}

A key result of Helfgott \cite{He1} shows that generating sets of
$SL(2,p)$ grow rapidly under multiplication. His result on diameters
is an immediate consequence.

\begin{thm}[Helfgott]
  Let $G=SL(2,p)$ and $A$ a generating set of $G$. Let $\delta$ be a
  constant, $0<\delta<1$.
  \begin{enumerate}[a)]
  \item \label{item:1} Assume that $|A|< |G|^{1-\delta}$. Then
    $$
    |A^3|>> |A|^{1+\e}
    $$
    where $\e$ and the implied constant depend only on $\delta$
  \item \label{item:2} Assume that $|A| > |G|^{1-\delta}$. Then
    $A^k=G$ where $k$ depends only on $\delta$.
  \end{enumerate}
\end{thm}

It was observed in \cite{NP} that a result of Gowers \cite{Gow}
implies that \ref{item:2}) holds for an arbitrary simple group of Lie
type $L$ with $k=3$ for some $\delta(r)$ which depends only on the Lie
rank $r$ of $L$ (see \cite{BNP} for a more detailed discussion).
Hence to complete the proof of our theorem on diameters it remains to
prove an analogue of the (rather more difficult) part \ref{item:1}) as
was done by Helfgott for the groups $SL(3,p)$ in \cite{He2}.

We prove the following.

\begin{thm} \label{simple-L}
  Let $L$ be a simple group of Lie type of rank $r$ and $A$ a
  generating set of $L$.
  Then either $A^3=L$ or
  $$
  |A^3|>> |A|^{1+\e}
  $$
  where $\e$ and the implied constant depend only on $r$.
\end{thm}

In fact, instead of concentrating on simple groups,
we work in the framework of arbitrary linear algebraic groups.
A version of the above theorem, valid for finite groups obtained from
an arbitrary reductive group, produces growth within certain normal
subgroups.

We also give some examples which show that in the above result the
dependence of $\e$ on $r$ is necessary. In particular we construct
generating sets of $SL(d,3)$ of size $1.1^d$ with $|A^3|< 10|A|$
for $d$ large enough.

The proofs of Helfgott combine group theoretic arguments with some
algebraic geometry, Lie theory and tools from additive
combinatorics such as the sum-product theorem of Bourgain, Katz,
Tao \cite{BKT} . Our argument relies on a deeper understanding of the
algebraic group theory behind his proofs and an extra trick,
but not on additive combinatorics.
Rather, another version of the above theorem gives growth results in
certain solvable groups, which can be interpreted as sum-product
type theorems. Such results also appear in \cite{He2} with
combinatorial proofs.

\section{A corollary on expanders}
   
Helfgott's work has been the starting point and inspiration of much
recent work by Bourgain, Gamburd, Sarnak and others.  Let
$S={g_1...g_k}$ be a subset of $SL(d,\BZ)$ and $\Lambda=\langle
S\rangle$ the subgroup generated by $S$.\break Assume that $\Lambda$
is Zariski dense in $SL(d)$.  According to the theorem of
Matthews-Vaserstein-Weisfeiler \cite{MVW} there is some integer $m_0$
such that $\pi_m(\Lambda)=SL(d,\BZ/m\BZ)$ assuming $(m,m_0)=1$.
Here $\pi_m$
denotes reduction $\mod m$. It was conjectured in \cite{Lu},
\cite{BGS} that the Cayley graphs 
$\Gamma\big(SL(d,\BZ/m\BZ),\pi_m(S)\big)$
form an
expander family, with expansion constant bounded below by a constant
$c=c(S)$.  This was verified in \cite{BG1} \cite{BG2} \cite{BGS} in
many cases when $d=2$ and in \cite{BG3} for $d>2$ and moduli of the
form $p^n$ where $n\to\infty$ and $p$ is a sufficiently large prime.

In \cite{BG3} Bourgain and Gamburd also prove the following

\begin{thm}[Bourgain, Gamburd]
  Assume that the analogue of Helfgott's theorem on growth holds for
  $SL(d,p)$, $p$ a prime. Let $S$ be a finite subset of $SL(d,Z)$
  generating a subgroup $\Lambda$ which is Zariski dense in $SL(d)$.
  Then the family of Cayley graphs $\Gamma(SL(d,p),\pi_p(S))$ forms an
  expander family as $p\to\infty$. The expansion coefficients are
  bounded below by a positive number $c(S)> 0$.
\end{thm}

By our result on growth the condition of the above theorem is
satisfied hence the conjecture is proved for prime moduli.

For $d=2$ Bourgain, Gamburd and Sarnak \cite{BGS} proved that the
conjecture holds for square free moduli. This result was used in
\cite{BGS} as a building block in a combinatorial sieve method
for primes and almost primes on orbits of various subgroups of
$GL(2,\BZ)$ as they act on $\BZ^n$.

Very recently P. Varj\'u \cite{Va} has shown that if the analogue of
Helfgott's theorem holds for $SL(d,p)$, $p$ prime then the conjecture
holds for square free moduli.  Hence our results constitute a major
step towards obtaining a generalisation to Zariski dense subgroups of
$SL(d,\BZ)$ and to other arithmetic groups.

\section{Methods}

% \begin{thm}
%   Let $L\le GL(n)$ be a reductive algebraic group and 
%   $N\triangleleft L$ a closed normal subgroup, both defined over the
%   finite field \Fq.
%   We denote by $L(\Fq)$ the finite group of
%   \Fq-points in $L$.
%   Let $A$ be a generating set of $L(\Fq)$ and $\delta$ a constant with
%   $0<\delta<1$. Assume that $|A|<\big|L(\Fq)\big|^{1-\delta}$.
%   Then there is a closed subgroup $M\le N$ which is normal in $L$
%   such that
%   $$
%   \big|A^k\cap M\big|>>|A\cap N|^{(1+\e)\dim(M)/\dim(N)}
%   $$
%   where $k$, $\e$ and the implied constant depends only on
%   $\delta$, $\dim(L)$ and the degrees of $L$ and $N$ (as subvarieties).
% \end{thm}

We prove various results which say that if $L$ is a ``nice''
subgroup of an algebraic group $G$ generated by a set $A$ then $A$
grows in some sense. In particular, if $G$ is reductive then we have
the following:

\begin{thm}\label{reductive-L}
  Let $G,H\le GL(d,\Fqc)$ be reductive algebraic groups defined over the
  finite field \Fq
  such that $H$ normalises $G$ and the centraliser $\CC_G(H)$ is finite. 
  There are constants $m$, $\e$ depending on $\dim(G)$,
  and a constant $K$ depending on  $\delta$, $\dim(G)$
  and the degree of $G$ (as a subvariety in $GL(d)$) such that
  whenever  $A$ is a generating set of $H(\Fq)$, the finite group of
  \Fq-points of $H$, with 
  $K\le|A|$
  then either the commutator $\big[G(\Fq),G(\Fq)\big]$
  is contained in $A^3$
  or there is a closed subgroup $N\le G$ normalised by $H$ such that
  $$
  \big|A^m\cap N\big| >>
  \big|A\cap G\big|^{(1+\e)\dim(N)/\dim(G)} \;.
  $$  
\end{thm}

In the case when $G=H$ is simple, one has $N=G$ necessarily,
hence Theorem~\ref{simple-L} follows for the Chevalley groups.
One can modify Theorem~\ref{reductive-L} so that it applies to all
finite simple groups of Lie type,
and one can easily generalise it to a large number of non-reductive
groups.

Let us outline the proof in the simplest case, when $A$ generates
$L=SL(d,q)$, $q$ a prime-power.  Assume that ``$A$ does not grow''
i.e. $|AAA|$ is not much larger than $|A|$.  Using an "escape from
subvarieties" argument is shown in \cite{He2} that if $T$ is a maximal
torus in $L$ then $|T \cap A|$ is not much larger than $|A|^{1/(d+1)}$ .
This is natural to expect for dimensional reasons since
$\dim(T)/\dim(L)=(d-1)/(d^2-1)=1/(d+1)$

This $T$ is equal to $L\cap\bar T$ where $\bar T$ is a maximal torus
of the algebraic group $SL(d,\Fqc)$.
Let $T_r$ denote the set of regular semisimple elements in $T$.  Note
that $T\setminus T_r$ is contained in a subvariety $V\subsetneq\bar T$
of dimension $d-2$.
We use a rather more general escape argument
to show that $|T\setminus T_r \cap A|$ is not much larger than
$$
|A|^{\dim(V)/\dim(L)} = |A|^{1/(d+1) - 1/(d^2 -1)} \;.
$$

By \cite{He2} or by our escape argument $A$ does contain regular
semisimple elements.  If $a$ is such an element then 
consider the map $SL(d)\to SL(d)$, $g\to g^{-1}ag$.
The image of this map is contained in a subvariety of dimension
$d^2-1-(d-1)$ since $\dim\big(\CC_{SL(d)}(a)\big)=d-1$.
By the escape argument
we obtain that for the conjugacy class $\cl(a)$ of $a$ in $L$,
$\big|\cl(a) \cap A^{-1}aA\big|$ is not much larger than
$|A|^{(d^2-d)/(d^2-1)}$. Now $\big|\cl(a) \cap A^{-1}aA\big|$ is at least
the number of cosets of the centraliser $C_L(a)$
which contain elements of $A$ . It follows
that $\big|AA^{-1} \cap C_L(a)\big|$ is not much smaller than $|A|^{1/(d+1)}$.
Of course $C_L(a)$ is just the (unique) maximal torus containing $a$.

Let us say that $A$ covers a maximal torus $T$ if $\big|T \cap A\big|$
contains a regular semisimple element.  We obtain the following
fundamental dichotomy.

\begin{prop}
  Assume that a generating set A does not grow
  \begin{enumerate}[i)]
  \item If $A$ does not cover a maximal torus $T$ then $\big|T \cap A\big|$ is
    not much larger than $|A|^{1/(d+1) -1/(d^2-1)}$.
  \item If $A$ covers $T$ then $\big|T \cap AA^{-1}\big|$ is not much smaller
    than $|A|^{1/(d+1)}$. In this latter case in fact 
    $\big|T_r\cap AA^{-1}\big|$ is not much smaller than $|A|^{1/(d+1)}$.
  \end{enumerate}
\end{prop}

It is well known that if $A$ doesn't grow then $B=AA^{-1}$ doesn't
grow either hence the Proposition applies to B.

Let us first assume that $B$ covers a maximal torus $T$ but does not
cover a conjugate $T'=g^{-1}Tg$ of $T$ for some element $g$ of $L$ .
Since $A$ generates $L$ we have such a pair of conjugate tori where
$g$ is in fact an element of $A$. Consider those cosets of $T'$ which
intersect $A$. Each of
the, say, $t$ cosets contains at most $|B \cap T'|$ elements of $A$
i.e. not much more than $|B|^{1/(d+1) -1/(d^2-1)}$ which in turn is
not much more than $|A|^{1/(d+1)-1/(d^2-1)}$.  Therefore $|A|$ is not
much larger than $t|A|^{1/(d+1) -1/(d^2 -1)}$.

On the other hand $A\big(A^{-1}(BB^{-1})A\big)$ has at least 
$t\big|T\cap BB^{-1}\big|$
elements which is not much smaller than $t|A|^{1/(d+1)}$.
Therefore $A\big(A^{-1}(AA^{-2}A)A\big)$ is not much smaller than
$|A|^{1+1/(d^2-1)}$ which contradicts the assumption that $A$ does not
grow.  We obtain that $B$ covers all conjugates of some maximal torus
$T$. Now the conjugates of the set $T_r$ are pairwise
disjoint (e.g. since two regular semisimple elements commute exactly
if they are in the same maximal torus).  The number of this tori is
$|L:N_L(T)| > c(d)|L:T|$ for some constant which depends only on $d$.
Each of them contains not much less than $|B|^{1/(d+1)}$ regular
semisimple elements of $BB^{-1}$. Altogether we see $|A|$ is not much
smaller than $q^{d^2-d}|A|^{1/(d+1)}$ and finally that $|A|$ is not
much less than $|L|$.  in this case by \cite{NP} we have $AAA=L$.

\end{document}